\documentclass[11pt,amssym,twoside]{article}
\usepackage{amssymb}
\newtheorem{thm}{Theorem}[section]

\newtheorem{conj}[thm]{Conjecture}
\newtheorem{prop}[thm]{Proposition}

\newtheorem{defn}[thm]{Definition}

\newtheorem{rem}[thm]{Remark}

\newcommand{\To}{\longrightarrow}
\newcommand{\Proj}{\mathbb P}
\newcommand{\Z}{\mathbb Z}

\newcommand{\Q}{\mathbb Q}
\begin{document}

\title{Arithmetic Teichmuller theory}
\author{Arash Rastegar}
%\address{}%
%\email{rastegar@sharif.edu}%

%\thanks{}%
%\subjclass{}%
%\keywords{}%

%\date{}%
%\dedicatory{}
%\commby{}%
% ----------------------------------------------------------------
\maketitle
\begin{abstract}
By Grothendieck's anabelian conjectures, Galois representations
landing in outer automorphism group of the algebraic fundamental
group which are associated to hyperbolic smooth curves defined
over number-fields encode all the arithmetic information of these
curves. The Goal of this paper is to develop an arithmetic
Teichmuller theory, by which we mean, introducing arithmetic
objects summarizing the arithmetic information coming from all
curves of the same topological type defined over number-fields.
We also introduce Hecke-Teichmuller Lie algebra which plays the
role of Hecke algebra in the anabelian framework.

\end{abstract}
%------------------------------------------------------------------
\section*{Introduction}
One canonically associates to a proper smooth curve $X$ which is geometrically reduced and is defined over a
number field $K$ a continuous group homomorphism
$$
\rho_X:Gal(\bar {K}/K) \To Out(\pi_1^{alg}(X))
$$
where $Out(\pi_1^{alg}(X)) $ denotes the quotient of the automorphism
group $Aut(\pi_1^{alg}(X)) $ by inner automorphisms of the algebraic
fundamental group. By a conjecture of Voevodski and Matsumoto the
outer Galois representation is injective when topological
fundamental group of $X$ is nonabelian. Special cases of this
conjecture are proved by Belyi for $\mathbb {P}^1-\{0,1,\infty\}$
[Bel], by Voevodski in cases of genus zero and one [Voe], and by
Matsumoto for affine $X$ using Galois action on profinite braid
groups [Mat].

By Grothendieck's anabelian conjectures, in the case of
hyperbolic curves, the Galois module structure of $Out(\pi_1(X))$
should inherit all the arithmetic information of the curve. In
particular, it should produce all points defined over
number-fields and should characterize the isomorphism class of
$X$ over $K$. The former is called Grothendieck's ``section
conjecture" and the latter is implied by Gorthendieck's ``Hom
conjecture" which is proved by Mochizuki [Moc]. Thus, given $X$
and $X'$ hyperbolic curves, the natural map
$$
Isom_{K}(X,X')\To Out_{Gal(\bar
{K}/K)}(Out(\pi_1^{alg}(X)),Out(\pi_1^{alg}(X')))
$$
is a one-to-one correspondence. Here $Out_{Gal(\bar {K}/K)}$
denotes the set of Galois equivariant isomorphisms between the two
profinite groups.

In this paper, we introduce an arithmetic structure summarizing
all such arithmetic information for hyperbolic smooth curves of
given topological type defined over $K$. More precisely, we shall
summarize all $\rho_X$ in a single Galois representation. This
would be the beginning of arithmetic Teichmuller theory.

From now on, we assume that $2g-2+n>0$ to ensure hyperbolicity.
The stack $M_{g,n}$ is defined as the moduli stack of $n$-pointed
genus $g$ curves. By a family of $n$-pointed genus $g$ curves
over a scheme $S$, we mean a proper smooth morphism $C\to S$
whose fibers are proper smooth curves of genus $g$, together with
$n$ sections $s_i:S\to C$ for $i=1,...,n$ whose images do not
intersect.

The moduli stack $M_{g,n}$ is an algebraic stack over $Spec(\Z)$.
One can define the etale fundamental group of the stack $M_{g,n}$
in the same manner one defines etale fundamental group of
schemes. Oda showed that the etale homotopy type of the algebraic
stack $M_{g,n}\otimes \bar {\Q}$ is the same as the analytic
stack $M^{an}_{g,n}$ and its algebraic fundamental group is
isomorphic to the completion $\widehat {\Gamma_{g,n}}$ of the
Teichmuller modular group, or the mapping class group of
$n$-punctured genus $g$ Riemann surfaces [Oda]:
$$
\pi_1^{alg}(M_{g,n}\otimes \bar {\Q}) \cong \widehat
{\Gamma_{g,n}}.
$$
Triviality of $\pi_2$ implies exactness of the following short
sequence for the universal family $C_{g,n}\To M_{g,n}$ over the
moduli stack
$$
0\To \pi_1^{alg}(X,b)\To \pi_1^{alg}(C_{g,n},b)\To
\pi_1^{alg}(M_{g,n},a) \To 0
$$
where $X$ is the fiber on $a$ and $b$ is a point on $C_{g,n}$.
Using this exact sequence, one defines the arithmetic universal
monodromy representation
$$
\rho_{g,n}:\pi_1^{alg}(M_{g,n},a) \To Out(\pi_1^{alg}(X))
$$
In fact, after restriction to  $\pi_1^{alg}(M_{g,n}\otimes \bar
{\Q},a)$, this is the completion of the natural map
$$
\Gamma_{g,n}\To Out(\Pi_{g,n})
$$
where $\Pi_{g,n}$ denotes the topological fundamental group of the
curve of genus $g$ with $n$ punctures.

On can think of the Galois module structure of
$Out(\pi_1^{alg}(M_{g,n}\otimes \bar{\Q}))$ as a replacement for
the Teichmuller space. This object has the information of all
outer representations associated to smooth curves over $\Q$.
Indeed, by fixing such a hyperbolic curve $X$ of genus $g$ with
$n$ punctures defined over $\Q$, we have introduced a rational
point on the moduli stack $a\in M_{g,n}$ and thus a Galois
representation
$$
Gal(\bar{\Q}/\Q)\To \pi_1^{alg}(M_{g,n},a)
$$
which splits the following short exact sequence
$$
0\To \pi_1^{alg}(M_{g,n}\otimes \bar {\Q},a)\To
\pi_1^{alg}(M_{g,n},a) \To Gal(\bar {\Q}/\Q)\To 0.
$$
Composing with the arithmetic universal monodromy representation,
we get
$$
Gal(\bar{\Q}/\Q)\To Out(\pi_1^{alg}(X))
$$
which recovers the canonical outer representation associated to
$X$. Therefore, the following universal Galois representation
$$
\rho_{univ}:Gal(\bar{\Q}/\Q)\To Out(\pi_1^{alg}(M_{g,n}\otimes
\bar{\Q}))=Out(\widehat {\Gamma_{g,n}})
$$
is the arithmetic analogue of the Teichmuller space. Here, we
have assumed that, one can treat $M_{g,n}$ as an anabelian space
whose arithmetic is governed by Grothendieck's anabelian
conjectures, as was expected by Grothendieck.

Having this picture in mind, we try to translate this arithmetic
information to the language of Lie algebras in order to make it
more accessible computationally.

%------------------------------------------------------------------
\section{Background material}

The study of outer representations of the Galois group has two
origins. One root is the theme of anabelian geometry introduced by
Grothendieck [Gro] which lead to results of Nakamura, Tamagawa and
Mochizuki who solved the problem in dimension one [Moc]. The
second theme which is originated by Deligne and Ihara
independently deals with Lie-algebras associated to the pro-$l$
outer representation [Del] [Iha]. This lead to a partial proof of
a conjecture by Deligne [Hai-Mat]. In the first part, we will
review the weight filtration introduced by Oda (after Deligne and
Ihara) and a circle of related results.

% ----------------------------------------------------------------
\subsection{Weight filtration on $\widetilde {Out}(\pi^l_1(X))$}

By a fundamental result of Grothendieck [Gro] the pro-$l$
geometric fundamental group $\pi^l_1(X_{\bar {K}})$ of a smooth
algebraic curve $X$ over $\bar {K}$ is isomorphic to the pro-$l$
completion of its topological fundamental group, after extending
the base field to the field of complex numbers. The topological
fundamental group of a Riemann surface of genus $g$ with $n$
punctured points has the following standard presentation:
$$
\Pi_{g,n}\cong <a_1,...,a_g,b_1,...,b_g,c_1,...,c_n|\prod_{i=1}^g
[a_i,b_i]\prod_{j=1}^n c_j=1>.
$$

Let $Aut(\pi^l_1(X))$ denote the group of continuous
automorphisms of the pro-$l$ fundamental group of $X$ and
$Out(\pi^l_1(X))$ denote its quotient by the subgroup of inner
automorphisms. We will induce filtrations on particular subgroups
of these two groups.

Let $\bar X$ denote the compactification of $X$ obtained by
adding finitely many points. $\bar X$ is still defined over $K$.
Let $\lambda:Aut(\pi^l_1(\bar {X}))\to GL(2g,\Z_l)$ denote the map
induced by abelianization. The natural actions of
$Aut(\pi^l_1(\bar X))$ on cohomology groups $H^i(\pi^l_1(\bar
X),\Z_l)$ are compatible with the non-degenerate alternating form
defined by the cup product:
$$
H^1(\pi^l_1(\bar X),\Z_l)\times H^1(\pi^l_1(\bar X),\Z_l)\to
H^2(\pi^l_1(\bar X),\Z_l)\cong \Z_l.
$$
This shows that the image of $\lambda$ is contained in
$GSp(2g,\Z_l)$. One can prove that $\lambda$ is surjective and if
$\tilde {\lambda}$ denote the natural map
$$
\tilde{\lambda}:Out(\pi^l_1(\bar {X})) \longrightarrow
GSp(2g,\Z_l)
$$
there are explicit examples showing that the Galois representation
$\rho^l_X\circ \tilde{\lambda}$ does not fully determine the
original anabelian Galois representation [Asa-Kon].

Let $\widetilde {Aut}(\pi^l_1(X))$ denote the Braid subgroup of
$Aut(\pi^l_1(X))$ which consists of those elements taking each
$c_i$ to a conjugate of a power $c_i^{\sigma}$ for some $\sigma$
in $\Z_l^*$. There is a natural surjective map
$$
\pi_X:\widetilde {Aut}(\pi^l_1(X)) \longrightarrow
Aut(\pi^l_1(\bar {X}))
$$
Oda uses $\tilde {\lambda}$ to define and study natural
filtrations on $\widetilde{Aut}(\pi^l_1(X))$ and
$\widetilde{Out}(\pi^l_1(X))$. In the special case of
$X=\Proj^1-\{0,1,\infty\}$ this is the same filtration as the
filtration introduced by Deligne and Ihara. This filtration is
also used by Nakamura in bounding Galois centralizers [Nak].
Consider the central series of the pro-$l$ fundamental group
$$
\pi^l_1(X)=I^1\pi^l_1(X)\supset I^2\pi^l_1(X) \supset ... \supset
I^m\pi^l_1(X) \supset ...
$$
and let $I^1 Aut(\pi^l_1(X))$ denote the kernel of $\pi_X\circ
\lambda$. The central series filtration is not the most
appropriate for non-compact $X$. In general, we consider the
weight filtration, namely the fastest decreasing central
filtration such that
$$
I^2
\pi^l_1(X)=<[\pi^l_1(X),\pi^l_1(X)],c_1,...,c_n>_{\textrm{norm}}
$$
where $<.>_{\textrm{norm}}$ means the closed normal subgroup
generated by these elements. For $m\geq 3$ we define
$$
I^m \pi^l_1(X)=<[I^i \pi^l_1(X),I^j
\pi^l_1(X)]|i+j=m>_{\textrm{norm}}.
$$
The weight filtration induces a filtration on the automorphism
group of braid type by normal subgroups
$$
\widetilde {Aut}(\pi^l_1(X))=I^0 \widetilde
{Aut}(\pi^l_1(X))\supset I^1 \widetilde {Aut}(\pi^l_1(X)) \supset
...\supset I^m \widetilde {Aut}(\pi^l_1(X))\supset ...
$$
$$
I^m \widetilde {Aut}(\pi^l_1(X))=\{ \sigma \in \widetilde
{Aut}(\pi^l_1(X))|x^{\sigma}x^{-1}\in I^{m+1}\pi^l_1(X)\textrm{
for all } x\in\pi^l_1(X)\}
$$
\begin{prop} The weight filtration on $\widetilde
{Aut}(\pi^l_1(X))$ satisfies
$$
[I^m \widetilde {Aut}(\pi^l_1(X)),I^n \widetilde
{Aut}(\pi^l_1(X))]\subset I^{m+n}\widetilde {Aut}(\pi^l_1(X))
$$
for all $m$ and $n$, and induces a Lie-algebra structure on the
associated graded object $Gr^{\bullet}_I \widetilde
{Aut}(\pi^l_1(X))=\bigoplus gr^m \widetilde {Aut}(\pi^l_1(X))$.
The graded pieces
$$
gr^m \widetilde {Aut}(\pi^l_1(X))=I^m \widetilde
{Aut}(\pi^l_1(X))/I^{m+1}\widetilde {Aut}(\pi^l_1(X))
$$
are free $\Z_l$-modules of finite rank for all positive $m$.
\end{prop}
The weight filtration on the automorphism group of pro-$l$
fundamental group induces a filtration on the outer automorphism
group of braid type
$$
\widetilde {Out}(\pi^l_1(X))=I^0 \widetilde
{Out}(\pi^l_1(X))\supset I^1 \widetilde {Out}(\pi^l_1(X)) \supset
...\supset I^m \widetilde {Out}(\pi^l_1(X))\supset ...
$$
\begin{prop} The induced filtration on $\widetilde
{Out}(\pi^l_1(X))$ satisfies
$$
[I^m \widetilde {Out}(\pi^l_1(X)),I^n \widetilde
{Out}(\pi^l_1(X))]\subset I^{m+n}\widetilde {Out}(\pi^l_1(X))
$$
for all $m$ and $n$, and induces a Lie-algebra structure on the
associated graded object $Gr_I^{\bullet}\widetilde
{Out}(\pi^l_1(X))=\bigoplus gr^m \widetilde {Out}(\pi^l_1(X))$.
The graded pieces
$$
gr^m \widetilde {Out}(\pi^l_1(X))=I^m \widetilde
{Out}(\pi^l_1(X))/I^{m+1}\widetilde {Out}(\pi^l_1(X))
$$
are finitely generated $\Z_l$-module for all positive $m$.
\end{prop}
These two propositions are proved in [kon]. One can induce
filtrations on $Aut(\pi^l_1(X))$ and $Out(\pi^l_1(X))$ using the
natural surjection
$$
Aut(\pi^l_1(X))\longrightarrow GSp(2g,\Z_l)
$$
which is induced by the action of $Aut(\pi^l_1(X))$ on $gr^1
\pi^l_1(X)\cong \Z_l^{2g}.$

%---------------------------------------------------------------
\subsection{Graded pieces of $\widetilde{Out}(\pi^l_1(X))$}

The explicit presentation of the fundamental group of a Riemann
surface given in the previous section implies that $\pi^l_1(X)$ is
a one-relator pro-$l$ group and therefore very close to a free
pro-$l$ group. The groups $Aut(\pi^l_1(X))$ and $Out(\pi^l_1(X))$
also look very similar to automorphism group and outer
automorphism group of a free pro-$l$ group [Rib-Zal]. This can be
shown more precisely in the particular case of
$\widetilde{Out}(\pi^l_1(X))$.

The graded pieces of $Gr^{\bullet}_I \widetilde{Out}(\pi^l_1(X))$
can be completely determined in terms of the graded pieces of
$Gr^{\bullet}_I \pi^l_1(X)$ which are free $\Z_l$-modules. In
fact, $Gr^{\bullet}_I \pi^l_1(X)$ is a free Lie-algebra over
$\Z_l$ generated by images of $a_i$'s and $b_i$'s in $gr^2
\pi^l_1(X)$ for $1\leq i\leq g$ and $c_j$'s in $gr^3 \pi^l_1(X)$
for $1\leq j\leq n$. We denote these generators by
$\bar{a}_i$,$\bar{b}_i$ and $\bar{c}_j$ respectively.

Let $g_m$ denote the following injective $\Z_l$-linear
homomorphism
$$
g_m:gr^m \pi^l_1(X) \longrightarrow (gr^{m+1} \pi^l_1(X))^{2g}
\times (gr^m \pi^l_1(X))^n
$$
$$
g\mapsto ([g,\bar{a}_i])_{1\leq i\leq g}\times
([g,\bar{b}_i])_{1\leq i\leq g}\times (g)_{1\leq j\leq n}
$$
and $f_m$ denote the following surjective $\Z_l$-linear
homomorphism
$$
f_m:(gr^{m+1} \pi^l_1(X))^{2g} \times (gr^m \pi^l_1(X))^n
\longrightarrow gr^{m+2} \pi^l_1(X)
$$
$$
(r_i)_{1\leq i\leq g}\times (s_i)_{1\leq i\leq g}\times
(t_j)_{1\leq j\leq n}\mapsto \sum_{i=1}^g
([\bar{a}_i,s_i]+[r_i,\bar{b}_i])+\sum_{j=1}^n [t_j,\bar{c}_j].
$$
\begin{prop} The graded pieces of $Gr^{\bullet}_I
\widetilde{Out}(\pi^l_1(X))$ fit into the following short exact
sequence of $\Z_l$-modules
$$
gr^m\widetilde{Out}(\pi^l_1(X)) \hookrightarrow (gr^{m+1}
\pi^l_1(X))^{2g} \times (gr^m \pi^l_1(X))^n/gr^m \pi^l_1(X)
\twoheadrightarrow gr^{m+2} \pi^l_1(X)
$$
where embedding of $gr^m \pi^l_1(X)$ inside $(gr^{m+1}
\pi^l_1(X))^{2g} \times (gr^m \pi^l_1(X))^n$ is defined by $g_m$
and the final surjection is induced by $f_m$. \end{prop}

This computational tool helps to work with the graded pieces of
$Gr^{\bullet}_I \widetilde{Out}(\pi^l_1(X))$ as fluent as the
graded pieces of $Gr^{\bullet}_I \pi^l_1(X)$. In particular, it
enabled Koneko to prove the following profinite version of the
Dehn-Nielson theorem [Kon]:
\begin{thm}(Koneko)
Let $X$ be a smooth curve defined over a number-field $K$ and let
$Y$ denote an embedded curve in $X$ obtained by omitting finitely
many $K$-rational points. Then the natural map
$$
\widetilde{Out}(\pi^l_1(Y))\longrightarrow
\widetilde{Out}(\pi^l_1(X))
$$
is a surjection.
\end{thm}
This can be easily proved by diagram chasing between the
corresponding short exact sequences for $Y$ and $X$. The above
exact sequence first appeared in the work of Ihara [Iha] and then
generalized by Asada and Koneko [Asa-Kon] [Kon].

%---------------------------------------------------------------
\subsection{Filtrations on the Galois group}

If all of the points in the complement $\bar {X}-X$ are
$K$-rational, then the pro-$l$ outer representation of the Galois
group lands in the braid type outer automorphism group
$$
\tilde{\rho}^l_X:Gal(\bar {K}/K) \longrightarrow
\widetilde{Out}(\pi^l_1(X))
$$
and the weight filtration on the pro-$l$ outer automorphism group
induce a filtration on the absolute Galois group mapping to
$\widetilde{Out}(\pi^l_1(X))$ and also an injection between
associated Lie algebras over $\Z_l$ defined by each of these
filtrations
$$
Gr^{\bullet}_{X,l}Gal(\bar {K}/K)\hookrightarrow Gr_I^{\bullet}
\widetilde{Out}(\pi^l_1(X)).
$$

\begin{prop}
Let $X$ and $Y$ denote smooth curves over $K$ and let $\phi:X\to
Y$ denote a morphism also defined over $K$. Then $\phi$ induces a
commutative diagram of Lie algebras
$$
\begin{array}{ccc}
Gr^{\bullet}_{X,l}Gal(\bar {K}/K)& \hookrightarrow &
Gr_I^{\bullet} \widetilde{Out}(\pi^l_1(X)) \\
\downarrow \phi_{**}& & \downarrow \phi_* \\
Gr^{\bullet}_{Y,l}Gal(\bar {K}/K) & \hookrightarrow &
Gr_I^{\bullet} \widetilde{Out}(\pi^l_1(Y))
\end{array}
$$
If $\phi$ induces a surjection on topological fundamental groups,
then $\phi_*$ and $\phi_{**}$ will also be surjective.
\end{prop}
\textbf{Proof} The claim is true because $\phi_*$ is Galois
equivariant and graded pieces of $Gr_I^{\bullet}
\widetilde{Out}(\pi^l_1(X))$ can be represented in terms of exact
sequences on graded pieces of $Gr_I^{\bullet} \pi^l_1(X)$ [Kon].
$\square$
\begin{prop}
Let $X$ be an affine smooth curve $X$ over $K$ whose complement
has a $K$-rational point. Then there is a morphism
$$
Gr^{\bullet}_{X,l}Gal(\bar {K}/K)\to Gr^{\bullet}_{\Proj^1-\{
0,1,\infty\},l}Gal(\bar {K}/K).
$$
\end{prop}
\textbf{Proof} This is a consequence of theorem 3.1 in
[Mat].$\square$
\begin{prop} There exists
a finite set of primes $S$ such that we have an isomorphism
$$
Gr^{\bullet}_{X,l}Gal(\bar {K}/K) \cong
Gr^{\bullet}_{X,l}Gal(K^{un}_S/K)
$$
where $Gal(K^{un}_S/K)$ denotes the Galois group of the maximal
algebraic extension unramified outside $S$
\end{prop}
\textbf{Proof.} Indeed, Grothendieck proved that the
representation $\tilde{\rho}^l_X$ factors through
$Gal(K^{un}_S/K)$ for a finite set of primes $S$. $S$ can be
taken to be primes of bad reduction of $X$ and primes over $l$
[Gro]. This is also proved independently by Ihara in the special
case of $X=\Proj^1-\{ 0,1,\infty\}$ [Iha].$\square$

The importance of this result of Grothensieck is the fact that
$Gal(K^{un}_S/K)$ is a finitely generated profinite group
[Neu-Sch-Wei] and therefore, the moduli of its representations is
a scheme of finite type.

For a Lie algebra $L$ over $\Z_l$ let $Der(L)$ denote the set of
derivations, which are defined to be $\Z_l$-linear homomorphisms
$D:L\to L$ with
$$
D([u,v])=[D(u),v]+[u,D(v)]
$$
for all $u$ and $v$ in $L$, and let $Inn(L)$ denote the set of
inner derivations, which are defined to be derivations with
$D(u)=[u,v]$ for some fixed $v\in L$. Then we have the following
Lie algebra version of the outer representation of the Galois
group
$$
Gr^{\bullet}_{X,l}Gal(\bar {K}/K)\to Der(Gr^{\bullet}_I
\pi^l_1(X))/Inn(Gr^{\bullet}_I \pi^l_1(X))
$$
$$
\sigma\in gr^m Gal(\bar {K}/K) \mapsto (u\mapsto
\tilde{\sigma}(\tilde{u}).\tilde{u}^{-1})\textrm{ mod }
I^{m+n+1}\pi^l_1(X)
$$
where $u\in gr^n \pi^l_1(X)$  with $\tilde{u}\in I^n \pi^l_1(X)$
and $\tilde{\sigma}\in I^m Gal(\bar {K}/K)$ is a lift of $\sigma$.
In fact, for a free graded algebra, one can naturally associate a
grading on the the algebra of derivations. Let $L=\oplus_i L^i$
be a free graded Lie algebra and let $D$ denote the derivation
algebra of $L$. Then define
$$
D^i=\{ d\in D| D(L^j)\subset L^{i+j} \}.
$$
Then every element $d\in D$ is uniquely represented in the form
$d=\sum_i d^i$ with $d^i\in D^i$ such that for any $f\in L$ the
component $d^i f$ vanishes for almost all $i$. One can prove that
$$
[D^i,D^j]\subset D^{i+j}\textrm{ and }[D_1,D_2]^k=\sum_{i+j=k}
[D_1^i,D_2^j].
$$
One can mimic the same construction on the graded algebra
associated to $\pi^l_1(X)$ to get a graded algebra of derivations
[Tsu].

One shall notice that in case $X=\Proj^1-\{ 0,1,\infty\}$ the
group $\pi^l_1(X)$ is the pro-$l$ completion of a free group with
two generators. Ihara proves that the associated Lie algebra
$Gr^{\bullet}_I(\pi^l_1(X))$ is also free over two generators say
$x$ and $y$ [Iha]. Now for $f$ in the $m$-th piece of the grading
of the Lie algebra, there is a unique derivation $D_f\in Der
(Gr^{\bullet}_I(\pi^l_1(X)))$ which satisfies $D_f(x)=0$ and
$D_f(y)=[y,f]$. One can show that $D_f(y)$ is non-zero for
non-zero $m$ and that for any $\sigma\in
Gr_I^{\bullet}(\pi^l_1(X))$ there exists a unique $f\in
Gr^{\bullet}_I(\pi^l_1(X))$ with image of $\sigma$ being equal to
$D_f$ [Mat2]. Now it is enough to let $\sigma=\sigma_m$ the Soule
elements, to get a non-zero image $D_{f_{\sigma_m}}$.

\begin{conj} (Deligne) The graded Lie algebra
$(Gr^{\bullet}_{\Proj^1-\{ 0,1,\infty\},l}Gal(\bar
{\Q}/\Q))\otimes \Q_l$ is a free graded Lie algebra over $\Q_l$
which is generated by Soule elements and the Lie algebra
structure is induced from a Lie algebra over $\Z$ independent of
$l$. \end{conj}
\begin{rem}
It is reasonable to expect freeness to hold for $(Gr^{\bullet}_{X
,l}Gal(\bar {\Q}/\Q))\otimes \Q_l$.
\end{rem}
This implies that the above graded Lie algebra representation is
also injective. Ihara showed that Soule elements do generate
$Gr_{X,l}^{\bullet}(Gal(\bar {\Q}/\Q))\otimes \Q_l$ if one assumes
freeness of this Lie algebra [Iha]. Hain and Matsumoto proved the
same result without assuming any part of Deligne's conjecture
[Hai-Mat].
% ----------------------------------------------------------------
\section{The completion of mapping class group}

The mapping class group $MC(X)$ of $X$ is defined as the factor of
the group of homeomorphisms of $X$ as a Riemann surface by the
subgroup of elements isotopic to identity. The mapping class
group is isomorphic to the group of outer automorphisms of the
topological fundamental group
$$
MC(X)\cong Out(\pi^{top}_{1}(X)).
$$

There has been many efforts to introduce a finite presentation
for this group. The ones introduced by Birman in 1974 for the case
of genus 2 look particularly simple. The generators
$\sigma_{1},...,\sigma_{5}$ together with the following relations
generate $MC(X_2).$
\medskip
\\ 
$$
\begin{array}{ll}
\sigma_i\sigma_j=\sigma_j\sigma_i, & \mbox{for } |i-j|\geq 2, \,  1\leq i,j\leq 5 \\
\sigma_i\sigma_{i+1}\sigma_i=\sigma_{i+1}\sigma_i\sigma_{i+1}, & 1\leq i\leq 4       \\
(\sigma_1\sigma_2 ...\sigma_5)^6=1 & \\
(\sigma_1\sigma_2\sigma_3\sigma_4\sigma_5^2\sigma_4\sigma_3\sigma_2\sigma_1)^2=1 &  \\
\sigma_1\sigma_2\sigma_3\sigma_4\sigma_5^2\sigma_4\sigma_3\sigma_2\sigma_1 \bot \sigma_i, & 1\leq i\leq 5
\end{array}
$$
\medskip
\\
where $x \bot y$ means that $x$ and $y$ commute. For $g\geq 3$
genus, Dehn 1938 [Deh], Lickorish 1965 [Lik], Hatcher and Thurston
1980 [Hat-Thu], Wajnryb 1983 [Waj], gave presentations of the mapping
class group. In all these presentations, the number of generators
increases with $g$. However, Suzuki in 1977 showed that one can
manage with four generators [Suz]. We give automorphisms in
$Aut(\pi^{top}_{1}(X))$ whose image in
$MC(X)=Aut(\pi^{top}_{1}(X))/Inn(\pi^{top}_{1}(X))$ generate the
mapping class group.
\bigskip
\\
$
\begin{array}{c}
\alpha_0: \left( \begin{array}{c}a_1\rightarrow b_1^{-1},
a_j\rightarrow a_j,j\neq 1\\
b_1\rightarrow b_1^{-1}a_1b_1,b_j\rightarrow b_j,j\neq 1
\end{array} \right) 
\end{array}
\medskip \\
\begin{array}{c}
\alpha_1: \left( \begin{array}{c}a_i\rightarrow a_{i+1},1\leq
i\leq g-1,a_g\rightarrow a_1\\
b_i\rightarrow b_{i+1},1\leq i\leq g-1,b_g\rightarrow b_1
\end{array} \right)
\end{array}
\medskip \\
\begin{array}{c}
\alpha_2: \left( \begin{array}{c}a_i\rightarrow a_{i},1\leq
i\leq g\\
b_1\rightarrow b_1^{-1},b_j\rightarrow b_j, 2\leq j\leq g
\end{array} \right)
\end{array}
\medskip \\
\begin{array}{c}
\alpha_3: \left( \begin{array}{c} a_2\rightarrow
b_2a_2(b_1^{-1}a_1b_1)(a_2^{-1}b_2^{-1}a_2)\\
a_j\rightarrow a_{j},j\neq 2\\
b_1\rightarrow b_1(a_2^{-1}b_2^{-1}a_2)\\
b_j\rightarrow b_{j},2\leq j\leq g
\end{array} \right) \\
\end{array}
\bigskip
$
\\ 
We define the groups $MC^-(X)$ and $MC^+(X)$ as the images of
$L^-$ and $L^+$ under the canonical homomorphism
$Aut(\pi^{top}_{1}(X))\longrightarrow MC(X)$. The generators for
$MC^-(X)$ and $MC^+(X)$ are also introduced by Suzuki in 1977. The
automorphisms $\alpha_1,\alpha_2$, and $\alpha_3$ generate
$MC^-(X)$. The generators of $MC^+(X)$ are the following. Here
$s_i$ denotes the word $b_i^{-1}a_i^{-1}b_ia_i$ for $1\leq i\leq
g$.
\bigskip
\\
$
\begin{array}{c}
\alpha_4: \left( \begin{array}{c}a_1\rightarrow
b_1^{-1}a_1^{-1}b_1, a_j\rightarrow a_j,2\leq j\leq g\\
b_1\rightarrow b_1^{-1}s_1^{-1}, b_j\rightarrow b_j,2\leq j\leq g
\end{array} \right)
\medskip \\
\end{array}
$\\
$
\begin{array}{c}
\alpha_5: \left( \begin{array}{c} a_1\rightarrow
s_1^{-1}a_2s_1,a_2\rightarrow a_1, a_j\rightarrow a_j,3\leq j\leq g\\
b_1\rightarrow s_1^{-1}b_2s_1,b_2\rightarrow b_1, b_j\rightarrow
b_j,3\leq j \leq g
\end{array} \right)
\medskip \\
\end{array}
$\\
$
\begin{array}{c}
\alpha_6: \left( \begin{array}{c}a_i\rightarrow a_{i},1\leq
i\leq g\\
b_1\rightarrow a_1b_1a_2^{-1}s_2(b_1^{-1}a_1^{-1}b_1)\\
b_2\rightarrow b_2a_2(b_1^{-1}a_1^{-1}b_1)a_2^{-1}\\
b_j\rightarrow b_{j},3\leq j\leq g
\end{array} \right)
\medskip \\
\end{array}
\bigskip
$
\\
By studying the mapping class group, we have considered
generators of $Out(\pi^{top}_{1}(X))$. In order to understand the
algebraic geometric analogue $Out(\pi^{alg}_{1}(X))$ we shall
study $Aut(\pi^{alg}_{1}(X))$ in more detail.

It is well known that, for a profinite group $G$ which admits a
fundamental system of open neighborhoods of the identity
consisting of characteristic subgroups, there exists a
topological isomorphism,
$$
Aut(G)\cong \lim_{\longrightarrow} Aut(G/U)
$$
where $U$ runs over open characteristic subgroups of $G$. We have
an injection $\pi^{top}_{1}(X)\longrightarrow\pi^{alg}_{1}(X))$.
An element of $Aut(\pi^{top}_{1}(X))$ fixes every characteristic
open subgroup $U$ and induces a compatible system of elements in
$Aut(A/U)$ for different $U$ and therefore an element of $\lim
Aut(G/U)$. We have constructed an injection
$$
Aut(\pi^{top}_{1}(X))\longrightarrow Aut(\pi^{alg}_{1}(X)).
$$
Inner automorphisms of $\pi^{top}_{1}(X)$ induce inner
automorphisms of the completion $\pi^{alg}_{1}(X)$. Thus we get a
second injection,
$$
Out(\pi^{top}_{1}(X))\longrightarrow Out(\pi^{alg}_{1}(X)).
$$
If we prove that $Aut(\pi^{alg}_{1}(X))$ is the profinite
completion of $Aut(\pi^{top}_{1}(X))$, we have shown that
$Out(\pi^{alg}_{1}(X))$ is the completion of $Out(\pi^{top}_{1}(X))\cong MC(X)$ 
in the profinite topology . It is enough to
show that $\pi^{alg}_{1}(X)$ has a fundamental system of open
characteristic subgroups which are completions of open subgroups
of $\pi^{top}_{1}(X)$ in the profinite topology. We know that
every automorphism of $\pi^{top}_{1}(X)$ is induced by an
automorphism of $F(a_1,...,a_g,b_1,...,b_g)$. So it is enough to
show that every free group has a fundamental system of open
characteristic subgroups. But this is proved to be true for a
finitely generated free group. Therefore we have a representation
of Galois group landing on the profinite completion of $MC(X)$
$$
\rho: Gal(\overline{\mathbb{Q}}/\mathbb{Q})\rightarrow
Out(\pi^{alg}_{1}(X))\cong \widehat{MC(X)}.
$$

Uchida in 1976 [Uch] and Ikeda in 1977 [Ike] proved that every
automorphism of $Gal(\overline{\mathbb{Q}}/\mathbb{Q})$ is inner.
Therefore the equivalence class of Galois representations
$$
\rho: Gal(\overline{\mathbb{Q}}/\mathbb{Q})\rightarrow
Out(\pi^{alg}_{1}(X))
$$
has only one element. This is unlike the case of $p$-adic Galois
representations associated to elliptic curves.

Translating the Galois representation from the language of
automorphisms to the language of mapping class group, provides us
a chance to geometrically define invariants of the Galois
representation. For example, one shall be able to give a purely
geometric definition of the conductor of a Galois representation.

%-----------------------------------------------------------------
\section{Arithmetic Teichmuller theory}
Here, we shall fulfill our promise of developing an arithmetic
Teichmuller theory. The universal representation
$$
\rho_{univ}:Gal(\bar{\Q}/\Q)\To Out(\widehat {\Gamma_{g,n}})
$$
as explained in the introduction summarizes the arithmetic
information coming from all curves of the same topological type
defined over number-fields. The aim is to translate this to the
language of Lie algebras. Although we lose some information, Lie
algebras are much more flexible for computations rather than than
Galois representations landing in outer automorphism groups.

In the first part, we introduce Hecke-Teichmuller Lie algebra
which plays the role of Hecke algebra in the anabelian framework.
Then, in the second section, we bring in the story of elliptic
curves and modularity. Finally, Galois action on the
Hecke-Teichmuller Lie algebra is related to elliptic curves.
% -----------------------------------------------------------------
\subsection{Hecke-Teichmuller graded Lie algebras}

There is no general analogue for Hecke operators in the context
of Lie algebras constructed in this manner. What we need is an
analogue of Hecke algebra which contains all the information of
Galois outer representations associated to elliptic curves. In
fact, we will provide an algebra containing such information for
hyperbolic smooth curves of given topological type.

From now on, we assume that $2g-2+n>0$. Note that elliptic curves
punctures at the origin are included. By composition with the
natural projection to outer automorphism group of the $l$-adic
completion $\Pi^l_{g,n}$ we get a representation
$$
\rho_{g,n}:\pi_1^{alg}(M_{g,n},a) \To Out(\pi_1^l(X)).
$$
This map induces filtrations on $\pi_1^{alg}(M_{g,n})$ and its
subgroup $\pi_1^{alg}(M_{g,n}\otimes \bar{\Q})$ and an injection
of $\Z_l$-Lie algebras
$$
Gr^{\bullet}_{X,l} \pi_1^{alg}(M_{g,n}\otimes
\bar{\Q})\hookrightarrow Gr^{\bullet}_{X,l} \pi_1^{alg}(M_{g,n}).
$$
It is conjectured by Oda and proved by a series of papers by
Ihara, Matsumoto, Nakamura and Takao that the cokernel of the
above map after tensoring with $\Q_l$ is independent of $g$ and
$n$ [Iha-Nak] [Mat2] [Nak]. Note that $M_{0,3}\cong Spec(\Q)$.

\begin{defn}
We define the Hecke-Teichmuller Lie algebra to be the image of the
following morphism of graded Lie algebras
$$
Gr^{\bullet}_{X,l} \pi_1^{alg}(M_{g,n}) \To Gr^{\bullet}_I
\widetilde {Out}(\Pi^l_{g,n}).
$$
\end{defn}

Note that, the filtration induced on the Galois group by
$Out(\Pi^l_{g,n})$ coincides with the filtration coming from
$\pi_1^{alg}(M_{g,n},a)$ by the Galois representation associated
to the corresponding curve and we have morphisms
$$
Gr^{\bullet}_{X ,l}Gal(\bar {\Q}/\Q)\To Gr^{\bullet}_{X,l}
\pi_1^{alg}(M_{g,n})\To Gr^{\bullet}_I \widetilde
{Out}(\Pi^l_{g,n}).
$$
We expect Hecke-Teichmuller Lie algebra to serve the role of
Hecke algebra in proving modularity results for elliptic curves
or other motivic objects.
%-------------------------------------------------------------------
\subsection{Galois representations associated to elliptic curves}

The method of proving modularity results by finding isomorphisms
between Hecke algebras and universal deformation rings as
originated by Wiles [Wil], can be reformulated in the language of
Lie algebras. One can find a canonical graded representation of
the Galois group to Hecke-Teichmuller Lie algebra which contains
all the information of modular Galois representations.

Let us first reformulate the theory of Galois representations in
the language Lie algebras. We start with elliptic curves. To each
elliptic curve $E$ defined over $\Q$ which has a rational point,
one associates a Galois outer representation
$$
Gal(\bar {\Q}/\Q) \longrightarrow Out(\pi^l_1(E-\{ 0\})).
$$
By analogy to Shimura-Taniyama-Weil conjecture, we expect this
representation be encoded in the representations
$$
Gal(\bar {\Q}/\Q) \longrightarrow Out(\pi^l_1(Y_0(N)))
$$
associated to modular curves $Y_0(N)$ which have a model over
$\Q$. By $Y_0(N)$ we mean the non-compactified modular curve of
level $N$ which is given as the quotient of the upper half-plane
by the congruence subgroup $\Gamma_0(N)$ of $SL_2(\Z)$ consisting
of matrices which are upper triangular modulo $N$.

For any smooth curve $X$ defined over $\Q$ the outer automorphism
group of braid type acts on $Gr^{\bullet}_I \widetilde
{Out}(\pi^l_1(X))$ by conjugation and therefore for each $m$ we
get a Galois representation
$$
Gal(\bar {\Q}/\Q) \longrightarrow Aut(gr^m \widetilde
{Out}(\pi^l_1(X))).
$$
In grade zero, we recover the usual abelian Galois representation
and in higher grades on can canonically construct this
representation by the grade-zero standard representation. Indeed,
for each $m\geq 1$ the isomorphism in proposition 1.3 is
$\widetilde {Out}(\pi^l_1(X))$-equivariant. From this we can
determine the representation from the inner action of $\widetilde
{Out}(\pi^l_1(X))$ on $Gr^{\bullet}_I \pi^l_1(X)$. This action is
fully determined by the grade-zero action. Therefore, the Galois
representations
$$
Gal(\bar {\Q}/\Q) \longrightarrow Aut(gr^m
\widetilde{Out}(\pi^l_1(X)))
$$
are all determined by the abelian Galois representation
associated to $X$ over $\Q$. Together with Shimura-Tanyama-Weil
conjecture proved by Wiles and his collaborators [Wil],[Tay-Wil],
[Bre-Con-Dia-Tay] we get the following Lie algebra version of
Shimura-Taniyama-Weil conjecture
\begin{thm}
Let $E$ be an elliptic curve over $\Q$ together with a rational
point $0\in E$. Then, the Galois representation
$$
Gal(\bar {\Q}/\Q)\longrightarrow Aut(Gr^{\bullet}_I
\widetilde{Out}(\pi^l_1(E-\{0\})))
$$
appear as direct summand of the Galois representation
$$
Gal(\bar{\Q}/\Q) \longrightarrow Aut(Gr^{\bullet}_I
\widetilde{Out}(\pi^l_1(Y_0(N))))
$$
for some level $N$.
\end{thm}

%-------------------------------------------------------------------
\subsection{Galois actions on Hecke-Teichmuller algebra}

The arithmetic analogue of the Teichmuller space which is the
universal Galois representation landing in outer automorphism
group of the algebraic fundamental group of
$M_{g,n}\otimes\bar{\Q}$
$$
\rho_{univ}:Gal(\bar{\Q}/\Q)\To
Out(\pi_1^{alg}(M_{g,n}\otimes\bar{\Q}))=Out(\widehat
{\Gamma_{g,n}})
$$
induces a Galois action on the corresponding Lie algebra
$$
Gal(\bar{\Q}/\Q) \longrightarrow Aut(Gr^{\bullet}_{I,l}
Out(\pi_1^{alg}(M_{g,n}\otimes\bar{\Q})))
$$
and thus a Galois action on the Hecke-Teichmuller Lie algebra.
This is analogue of the Galois action on the abelian Hecke
algebra, which is used to prove modularity of Galois
representations in [Wil]. One can use the Hecke-Teichmuller Lie
algebra to prove that certain Galois actions on Lie algebras do
come from curves defined over number fields.

The general question would be recognition of Galois actions on Lie
algebras which are motivic or even, recognition of Lie
representations of Galois-Lie group which are motivic. In order
to answer this general question, one shall find appropriate
generalizations of arithmetic hyperbolicity and Grothendieck's
conjectures or even a motivic formulation of hyperbolicity in
higher dimensions.
%-------------------------------------------------------------------
\section*{Acknowledgements}
I would like to thank M. Kontsevich, O. Gabber, N. Nikolov, V.B.
Mehta, A.J. Parameswaran, D. Prasad, Y. Soibelman for enjoyable
conversations which led to my interest in outer representations. I
wish to thank Tata institute of fundamental research and
institute des hautes \`{e}tudes scientifiques for their warm
hospitality during which these conversations took place. I would
also like to thank Harish-Chandra research institute,
particularly R. Kulkarni and R. Ramakrishnan for warm hospitality
during the workshop on Teichmuller theory and moduli spaces, for
which this note is prepared.
%-------------------------------------------------------------------

Sharif University of Technology, e-mail: rastegar@sharif.edu


\begin{thebibliography}{999999}
%\bibliographystyle{amsplain}
%\bibliography{bibliography}

\bibitem[Asa-Kon]{Asa-Kon} M. Asada, M. Koneko; {\em On the automorphism
group of some pro-$l$ fundamental groups,} Advanced studies in
Pure Mathematics Vol. 12,137-159(1987).

\bibitem[Bre-Con-Dia-Tay]{Bre-Con-Dia-Tay} C. Breuil, B. Conrad, F. Diamond,
R. Taylor; {\em On the modularity of elliptic curves over $\Q$:
wilde 3-adic excercises,} J. Amer. Math. Society 14 no.
4,843-939(2001).

\bibitem[Deh]{Deh} M. Dehn; {\em Die Gruppen der Abbildungsklassen},
Acta Math. 69, 135-206. Zbl.19.253.

\bibitem[Del]{Del} P. Deligne; {\em Le groupe fundamental de la droite
projective moins trois points,} in {\rm Galois groups over $\Q$}
79-297, edited by Y. Ihara, J.P. Serre, Springer 1989.

\bibitem[Gro]{Gro} A. Grothendieck; {\em Revetement Etale et
Groupe Fondamental (SGA I),} LNM 224, Springer-Verlag 1971.

\bibitem[Hat-Thu]{Hat-Thu} Hatcher A., Thurston W.; {\em A presentation
for the mapping class group of a closed orientable surface},
Topology 19, 221-237. Zbl.447.57005.

\bibitem[Iha]{Iha} Y. Ihara; {\em Profinite braid groups, Galois
representations, and complex multiplications,} Ann. of Math.
123,43-106(1986).

\bibitem[Iha-Nak]{Iha-Nak} Y. Ihara, H. Nakamura; {\em On deformation
of mutually degenerate stable marked curves and Oda's problem,}
J. Reine Angenw. Math. 487,125-151 (1997)

\bibitem[Ik]{Ike} Ikeda M.; {\em Completeness of the absolute Galois
group of the rational number field}, J. Reine Angew. Math., 291,
1-22.

\bibitem[Lab]{Lab} J. Labute; {\em On the descending central
series of groups with single defining relation,} J. of Algebra,
14,16-23(1970).

\bibitem[Lik]{Lik} Likorish W.B.R.; {\em On the homeomorphisms of a
non-orientable surface}, Proc. Cam. Philos. Soc.
61, 61-65. Zbl.131.208.

\bibitem[Oda]{Oda} T. Oda; {\em Etale homotopy type of the moduli
spaces of algebraic curves,} in "Geometric Galois Actions I"
London Math. Soc. Lect. Note Ser. 242,85-95 (1997)

\bibitem[Mat]{Mat} M. Matsumoto; {\em On the Galois image in the
derivation algebra of $\pi_1$ of the projective line minus three
points,} Contemporary Math. 186,201-213(1995).

\bibitem[Mat2]{Mat2} M. Matsumoto; {\em Galois representations on
profinite braid groups on curves,} J. Reine Angenw. Math.
474,169-219 (1996)

\bibitem[Moc]{Moc} S. Mochizuki; {\em The profinite Grothendieck
conjecture for closed hyperbolic curves over number fields,} J.
Math. Sci. Univ. Tokyo 3,571-627(1996).

\bibitem[Nak]{Nak} H. Nakamura; {\em Coupling of universal monodromy
representations of Galois Techmuller modular groups,} Math. Ann
304,99-119 (1996)

\bibitem[Rib-Zal]{Rib-Zal} L. Ribes, P. Zalesskii; {\em Prifinite
Groups,} Springer 2000.

\bibitem[Tay-Wil]{Tay-Wil} R. Taylor, A. Wiles; {\em
Ring-theoretic properties of certain Hecke algebras,} Ann. of
Math. 142,553-572(1995).

\bibitem[Suz]{Suz} S. Suzuki;{\em On homeomorphisms of a
3-dimensional handle body}, Can. J. Math. 29, No.1, 111-124.
Zbl.339.57001.

\bibitem[Tsu]{Tsu} H. Tsungai; {\em On some derivations of Lie algebras
related to Galois representations,} Publ. RIMS. Kyoto Univ.
31,113-134(1995).

\bibitem[Uch]{Uch} K. Uchida; {\em Isomorphisms of Galois groups,} J.
Math. Soc. Japan, 28, 617-620.

\bibitem[Voe]{Voe} V.A. Voevodskii; {\em Galois representations
connected with hyperbolic curves,} Math. USSR Izvestiya
39,1281-1291 (1992)

\bibitem[Waj]{Waj} B. Wajnryb; {\em A simple presentation for the
mapping class group of an orientable surface}, Isr. J. Math.
45, No. 2, 3, 157-174. Zbl.533.57022.

\bibitem[Wil]{Wil} A. Wiles; {\em Modular elliptic curves and Fermat's last
theorem,} Ann. of Math. 142,443-551(1995).


\end{thebibliography}
\end{document}